\documentclass[a4paper,11pt]{amsart}
\usepackage{amssymb}

\numberwithin{equation}{section}
\newtheorem{theorem}{Theorem}[section]

\newtheorem{lem}[theorem]{Lemma}

\theoremstyle{remark}
\newtheorem{rem}[theorem]{Remark}

\author[J.~Benameur]{Jamel Benameur}
\address{Department of Mathematics, College of Science, King Saud University\\
Riyadh 11451, Kingdom of Saudi Arabia}
\email{\sl jbenameur@ksu.edu.sa}
\author[L.~Jlali]{Lotfi Jlali}
\address{Department of Mathematics, College of Science, King Saud University\\
Riyadh 11451, Kingdom of Saudi Arabia}
\email{\sl ljlali@ksu.edu.sa}

\title[On the blow up criterion of 3D-NSE in Sobolev-Gevrey spaces]
{On the blow up criterion of 3D-NSE in Sobolev-Gevrey spaces}
\date{\today}

\begin{document}
\begin{abstract}In \cite{JB1}, Benameur proved a blow-up result of the non regular solution of $(NSE)$ in
the Sobolev-Gevrey spaces. In this paper we improve this result, precisely we give an exponential
type explosion in Sobolev-Gevrey spaces with less regularity on the initial condition. Fourier analysis is used.
\end{abstract}



\maketitle
\tableofcontents


\section{Introduction}
The $3D$ incompressible Navier-Stokes equations are given by:

$$
\left\{
  \begin{array}{lll}
     \partial_t u
 -\nu\Delta u+ u.\nabla u  & =&\;\;-\nabla p\hbox{ in } \mathbb R^+\times \mathbb R^3\\
     {\rm div}\, u &=& 0 \hbox{ in } \mathbb R^+\times \mathbb R^3\\
    u(0,x) &=&u^0(x) \;\;\hbox{ in }\mathbb R^3,
  \end{array}
\right.
\leqno{(NSE)}
$$
where $\nu$ is the viscosity of fluid, $u=u(t,x)=(u_1,u_2,u_3)$ and $p=p(t,x)$ denote respectively the unknown velocity and the unknown pressure of the fluid at the point $(t,x)\in \mathbb R^+\times \mathbb R^3$, and $(u.\nabla u):=u_1\partial_1 u+u_2\partial_2 u+u_3\partial_3u$, while $u^0=(u_1^o(x),u_2^o(x),u_3^o(x))$ is a given initial velocity. If $u^0$ is quite regular, the divergence free condition determines the pressure $p$.\\

Our problem is the study of explosions for non-smooth solutions of $(NSE)$. In the literature, several authors studied this problem (see \cite{BKM,JB01,JB02,Le2}), however all the obtained results do not exceed $(T^*-t)^{-\sigma}$, where $T^*$ denotes the maximal time of existence. Recently, in \cite{JB1}, the author gives a positif answer on the question: Is the type of explosion is due to the chosen space or the nonlinear part of $(NSE)$? For this, he used the Sobolev-Gevrey spaces which are defined as follows; for $a, s\geq0$ and $\sigma\geq 1$,
$$H^s_{a,\sigma}(\mathbb R^3)=\{f\in L^2(\mathbb R^3);\;e^{a|D|^{1/\sigma}}f\in H^s(\mathbb R^3)\}$$
is equipped with the norm
$$\|f\|_{H^s_{a,\sigma}}^2=\|e^{a|D|^{1/\sigma}}f\|_{H^s}$$
and its associated inner product
$$\langle f/g\rangle_{{H}^s_{a,\sigma}(\mathbb{R}^3)}=\langle e^{a|D|^{\frac{1}{\sigma}}}f/e^{a|D|^{\frac{1}{\sigma}}}g\rangle_{{H}^s(\mathbb{R}^3)}.$$
Precisely, for  $a>0$, $s>3/2$, $\sigma>1$ and $u^0\in \big(H^s_{a,\sigma}(\mathbb R^3)\big)^3$ such that ${\rm div}\,u^0=0$. He proved that, there is a unique time $T^*\in(0,\infty]$ and a unique solution $u\in\mathcal C([0,T^*),H^s_{a,\sigma}(\mathbb R^3))$ to $(NSE)$ system such that $u\not\in\mathcal C([0,T^*],H^s_{a,\sigma}(\mathbb R^3))$. Moreover, if $T^*<\infty$, then
\begin{equation}\label{princip-I1}
c(T^*-t)^{-s/3}\exp\Big(aC(T^*-t)^{-\frac{1}{3\sigma}}\Big)\leq  \|u(t)\|_{H^s_{a,\sigma}},\,\forall t\in[0,T^*),
\end{equation}
where $c=c(s,u^0,\sigma)>0\;\;\;and\;\;\;C=C(s,u^0,\sigma)>0.$ The choice of this spaces is due to the scaling property: If $u(t,x)$ is a solution to $(NSE)$ system with the initial data $u^0(x)$, then for any $\lambda>0$, $\lambda u(\lambda ^2t,\lambda x)$ is a solution to $(NSE)$ system with the initial data $\lambda u^0(\lambda x)$. Also, we recall that the classical used spaces ${\bf B}(\mathbb R^3)$ for $(NSE)$ system satisfy the fundamental condition:
$$u^0(x)\in {\bf B}(\mathbb R^3)\Longleftrightarrow \lambda u^0(\lambda x)\in {\bf B}(\mathbb R^3),\;\forall \lambda>0,$$
which is not valid in the Sobolev-Gevrey spaces $H^s_{a,\sigma}(\mathbb R^3)$.\\

Our work is intended to improve the result in \cite{JB1}. Before generalizing this result, we recall the energy estimate: If $u\in{\mathcal C}([0,T],H^s(\mathbb R^3))$, with $s>5/2$, is a solution to $(NSE)$, then
\begin{equation}\label{ee1}
\|u(t)\|_{L^2}^2+2\nu\int_0^t\|\nabla u(z)\|_{L^2}^2dz=\|u^0\|_{L^2}^2.
\end{equation}
Now we are ready to state our main result.
\begin{theorem}\label{theo1}
Let $a>0$ and $\sigma>1$. Let $u^0\in \left(H^1_{a,\sigma}(\mathbb{R}^3)\right)^3$ be such that ${\rm div}\, u^0=0$, then there is a unique $T^*\in(0,\infty]$ and unique $u\in{\mathcal C}([0,T^*),H^1_{a,\sigma}(\mathbb R^3))$ solution to $(NSE)$ system such that $u\notin{\mathcal C}([0,T^*],H^1_{a,\sigma}(\mathbb R^3))$. If $T^*<\infty$, then
\begin{eqnarray}\label{ineq99}
\frac{c_1}{(T^*-t)^{\frac{2\sigma_0+1}{3\sigma}+\frac{1}{3}}}\exp\Big[\frac{ac_2}{(T^*-t)^{1/3\sigma}}\Big]\leq \|u(t)\|_{H^1_{a,\sigma}},
\end{eqnarray}
where $c_1=c_1(u^0,a,\sigma)>0$, $c_2=c_1(u^0,\sigma)>0$ and $2\sigma_0$ is the integer part of $2\sigma$.
\end{theorem}
The remainder of this paper is organized in the following way: In section $2$, we give some notations and important preliminary results. Section $3$ is devoted to prove that $(NSE)$ is well posed in $H^1_{a,\sigma}(\mathbb R^3))$. In section $4$, we prove the exponential type explosion of non regular solution to $(NSE)$ system.
\section{Notations and preliminary results}
\subsection{Notations}
In this section, we collect some notations and definitions that will be used later.\\
$\bullet$ The Fourier transformation is normalized as
$$
\mathcal{F}(f)(\xi)=\widehat{f}(\xi)=\int_{\mathbb R^3}exp(-ix.\xi)f(x)dx,\,\,\,\xi=(\xi_1,\xi_2,\xi_3)\in\mathbb R^3.
$$
$\bullet$ The inverse Fourier formula is
$$
\mathcal{F}^{-1}(g)(x)=(2\pi)^{-3}\int_{\mathbb R^3}exp(-i\xi.x)g(\xi)d\xi,\,\,\,x=(x_1,x_2,x_3)\in\mathbb R^3.
$$
$\bullet$ For $ s\in\mathbb R $, $H^s(\mathbb R^3)$ denotes the usual non-homogeneous Sobolev space on $\mathbb R^3$ and  $\langle./.\rangle_{H^s(\mathbb R^3)}$ denotes the usual scalar product on $H^s(\mathbb R^3)$.\\
$\bullet$ For $ s\in\mathbb R $, $\dot{H}^s(\mathbb R^3)$ denotes the usual homogeneous Sobolev space on $\mathbb R^3$ and $\langle./.\rangle_{\dot{H}^s(\mathbb R^3)}$ denotes the usual scalar product on $\dot{H}^s(\mathbb R^3)$.\\
$\bullet$ The convolution product of a suitable pair of function $f$ and $g$ on $\mathbb R^3$ is given by
$$
(f\ast g)(x):=\int_{\mathbb R^3}f(y)g(x-y)dy.
$$
$\bullet$ If $f=(f_1,f_2,f_3)$ and $g=(g_1,g_2,g_3)$ are two vector fields, we set
$$
f\otimes g:=(g_1f,g_2f,g_3f),
$$
and
$$
{\rm div}\,(f\otimes g):=({\rm div}\,(g_1f),{\rm div}\,(g_2f),{\rm div}\,(g_3f)).
$$
$\bullet$ For $a>0$, $\sigma\geq1$, we denote the Sobolev-Gevrey space defined as follows
$$\dot H^1_{a,\sigma}(\mathbb R^3)=\{f\in{\mathcal S}'(\mathbb R^3);\;e^{a|D|^{1/\sigma}}f\in \dot H^1(\mathbb R^3)\}$$
which is equipped with the norm $$\|f\|_{\dot H^1_{a,\sigma}}=\|e^{a|D|^{1/\sigma}}f\|_{\dot H^1}$$
and the associated inner product  $$\langle f/g \rangle_{\dot H^1_{a,\sigma}}=\langle e^{a|D|^{1/\sigma}}f/e^{a|D|^{1/\sigma}}g \rangle_{\dot H^1}.$$
\subsection{Preliminary results}
In this section we recall some classical results and we give new technical lemmas.
\begin{lem}\label{lem1}(See \cite{HB})
Let $(s,t)\in{\mathbb{R}^2}$, such that $s<{\frac{3}{2}}$ and $s+t>0$. Then there exists a constant $C$, such that
$$
\|uv\|_{\dot{H}^{{s+t-{\frac{3}{2}}}}}\leq C(\|u\|_{\dot{H^s}}\|v\|_{\dot{H^t}}+\|u\|_{\dot{H^t}}\|v\|_{\dot{H^s}}).
$$
If $s<{\frac{3}{2}}$, $t<{\frac{3}{2}}$ and $s+t>0$,  then there exists a constant $C$, such that
$$
\|uv\|_{\dot{H}^{{s+t-{\frac{3}{2}}}}}\leq C \|u\|_{\dot{H^s}}\|v\|_{\dot{H^t}}.
$$
\end{lem}
\begin{lem}\label{lem2}(See \cite{JB1})
For $\delta>3/2$, we have
$$\|\widehat{f}\|_{L^1(\mathbb R^3)}\leq C_\delta\|f\|_{L^2}^{1-\frac{3}{2\delta}}\|f\|_{\dot H^\delta(\mathbb R^3)}^{\frac{3}{2\delta}},$$
with $C_\delta=2\sqrt{\frac{\pi}{3}}\Big((\frac{2\delta}{3}-1))^{3/(4\delta)}+
(\frac{2\delta}{3}-1)^{-1+\frac{3}{4\delta}}\Big).$\\
Moreover, for all $\delta_0>3/2$, there is $M(\delta_0)>0$ such that
$$C_\delta\leq M(\delta_0),\;\;\;\forall\delta\geq\delta_0.$$
\end{lem}
\begin{lem}\label{lem4}
For every $a>0$, $\sigma\geq1$, and for every $f, g \in\dot{H}^1_{a,\sigma}(\mathbb R^3)$, we have $fg\in\dot{H}^1_{a,\sigma}(\mathbb R^3)$ and
$$
\|fg\|_{\dot{H}^1_{a,\sigma}(\mathbb R^3))}\leq 16 \left(\|\mathcal{F}(e^{\frac{a}{\sigma}|D|^{\frac{1}{\sigma}}}f)\|_{ L^1}\|g\|_{\dot{H}^1_{a,\sigma}}+\|\mathcal{F}(e^{\frac{a}{\sigma}|D|^{\frac{1}{\sigma}}}g)\|_{ L^1}\|f\|_{\dot{H}^1_{a,\sigma}}\right)
$$
\end{lem}
\noindent{\it Proof of the lemma \ref{lem4}.}
We have
\begin{eqnarray*}
\|fg\|_{\dot{H}^1_{a,\sigma}(\mathbb R^3))}^2&=&\int_{\xi}|\xi|^2e^{2a|\xi|^{\frac{1}{\sigma}}}|\widehat{fg}(\xi)|^2\\&=&
\int_{\xi}|\xi|^2e^{2a|\xi|^{\frac{1}{\sigma}}}|\widehat{f}\ast\widehat{g}(\xi)|^2\\&\leq&\int_{\xi}|\xi|^2e^{2a|\xi|
^{\frac{1}{\sigma}}}(|\widehat{f}|\ast|\widehat{g}|)^2(\xi)\\&\leq&\int_{\xi}\left(\int_{\eta}|\xi|e^{a|\xi|^{\frac{1}{\sigma}}}
|\widehat{f}(\xi-\eta)||\widehat{g}(\eta)|\right)^2\\&\leq&\int_{\xi}\left(\int_{|\eta|<|\xi-\eta|}
|\xi|e^{a|\xi|^{\frac{1}{\sigma}}}|\widehat{f}(\xi-\eta)||\widehat{g}(\eta)|+\int_{|\eta|>|\xi-\eta|}|\xi|e^{a|\xi|^{\frac{1}{\sigma}}}
|\widehat{f}(\xi-\eta)||\widehat{g}(\eta)|\right)^2.
\end{eqnarray*}
By the elementary inequality
$$(1+b)^\theta\leq 1+\theta b^\theta,\;\forall b,\theta\in[0,1],$$
we can deduce that $$|\xi|^{\frac{1}{\sigma}}\leq \max(|\xi-\eta|,|\eta|)^{\frac{1}{\sigma}}+\frac{1}{\sigma}
\min(|\xi-\eta|,|\eta|)^{\frac{1}{\sigma}},\;\;\forall \xi,\eta\in\mathbb R^3.$$
Therefore
\begin{eqnarray*}
\|fg\|_{\dot{H}^1_{a,\sigma}(\mathbb R^3))}^2&\leq& 4\int_{\xi}(I_1+I_2)^2d\xi\\&\leq&16\int_{\xi}I_1^2(\xi)+I_2^2(\xi)d\xi,
\end{eqnarray*}
where
$$
I_1(\xi)=\int_{\eta}|\xi-\eta|e^{a|\xi-\eta|^{\frac{1}{\sigma}}}|\widehat{f}(\xi-\eta)|e^{\frac{a}{\sigma}
|\eta|^{\frac{1}{\sigma}}}|\widehat{g}(\eta)|
$$
$$
I_2(\xi)=\int_{\eta}e^{\frac{a}{\sigma}|\xi-\eta|^{\frac{1}{\sigma}}}|\widehat{f}(\xi-\eta)|e^{a|\eta|^{\frac{1}{\sigma}}}
|\eta||\widehat{g}(\eta)|.
$$
Using the following notations:
$$
\alpha_1(\xi)=|\xi|e^{a|\xi|^{\frac{1}{\sigma}}}|\widehat{f}(\xi)|
$$
$$
\alpha_2(\xi)=e^{\frac{a}{\sigma}|\xi|^{\frac{1}{\sigma}}}|\widehat{f}(\xi)|
$$
$$
\beta_1(\xi)=|\xi|e^{a|\xi|^{\frac{1}{\sigma}}}|\widehat{g}(\xi)|
$$
$$
\beta_2(\xi)=e^{\frac{a}{\sigma}|\xi|^{\frac{1}{\sigma}}}|\widehat{g}(\xi)|,
$$
we obtain
\begin{eqnarray*}
\|fg\|_{\dot{H}^1_{a,\sigma}(\mathbb R^3))}^2\leq16(\|\alpha_1\ast\beta_2\|_{L^2}^2+\|\alpha_2\ast\beta_1\|_{L^2}^2).
\end{eqnarray*}
Young's inequality yields
\begin{eqnarray*}
\|fg\|_{\dot{H}^1_{a,\sigma}(\mathbb R^3))}^2\leq16(\|\alpha_1\|_{L^2}^2\|\beta_2\|_{L^1}^2+\|\alpha_2\|_{L^1}^2\|\beta_1\|_{L^2}^2),
\end{eqnarray*}
and the proof of lemma \ref{lem4} is finished.\\
\begin{lem}\label{lem53} Let $a>0$, $\sigma\geq1$ and $s\geq1$, then there is a constant $c=c(s,a,\sigma)$ such that for all $f\in \dot H^1_{a,\sigma}(\mathbb R^3)$
\begin{equation}
\|f\|_{\dot H^s}\leq c\|f\|_{\dot H^1_{a,\sigma}}.
\end{equation}
\end{lem}
\noindent{\it Proof of the lemma \ref{lem53}.} Let $k_0\in\mathbb Z^+$ such that
$$\frac{k_0}{2\sigma}\leq s-1<\frac{k_0+1}{2\sigma}.$$
Using the fact that
$$|\xi|^{2s-2}\leq |\xi|^{\frac{k_0}{\sigma}}+|\xi|^{\frac{k_0+1}{\sigma}},\;\forall \xi\in\mathbb R^3,$$
we get
$$\|f\|_{\dot H^s}^2\leq \frac{(k_0+1)!}{(2a)^{k_0+1}}\big(2a+1\big)\int_{\mathbb R^3}\Big(\frac{(2a|\xi|^{1/\sigma})^{k_0}}{k_0!}+\frac{(2a|\xi|^{1/\sigma})^{k_0+1}}{(k_0+1)!}\Big)|\xi|^2|\widehat{f}(\xi)|^2\leq 2\frac{(k_0+1)!}{(2a)^{k_0}}\|f\|_{\dot H^1_{a,\sigma}}^2.$$
\begin{lem}\label{lem51}
For every $a>0$, $\sigma\geq1$, and for every $f, g \in H^1_{a,\sigma}(\mathbb R^3)$, we have $fg\in H^1_{a,\sigma}(\mathbb R^3)$ and
\begin{equation}\label{eq2a}
\|fg\|_{L^2(\mathbb R^3))}\leq C \left(\|f\|_{\dot{H}^1_{a,\sigma}}\|g\|_{L^2}+\|g\|_{\dot{H}^1_{a,\sigma}}\|f\|_{L^2}\right).
\end{equation}
\end{lem}
\noindent{\it Proof of the lemma \ref{lem51}.} We have
$$\|fg\|_{L^2(\mathbb R^3))}\leq C\Big(\|f\|_{L^\infty(\mathbb R^3))}\|g\|_{L^2(\mathbb R^3))}+\|f\|_{L^2(\mathbb R^3))}\|g\|_{L^\infty(\mathbb R^3))}\Big).$$
Using the classical interpolation inequality
$$\|h\|_{L^\infty}\leq c\|h\|_{\dot H^1}^{1/2}\|h\|_{\dot H^2}^{1/2}$$
and lemma \ref{lem53}, we can deduce inequality  (\ref{eq2a}).
\begin{rem} By lemmas \ref{lem4}-\ref{lem51}, we can deduce that $H^1_{a,\sigma}(\mathbb R^3)$ is an algebra for $a>0$ and $\sigma>1$.\\
\end{rem}
\begin{lem}\label{lem52} If $f\in H^1_{a,\sigma}(\mathbb R^3)$, then
\begin{equation}
\|f\|_{H^1_{a,\sigma}}^2\leq 2(e^{2a}+1)\Big(\|f\|_{L^2}^2+\|f\|_{\dot H^1_{a,\sigma}}^2\Big)\leq 4(e^{2a}+1)\|f\|_{H^1_{a,\sigma}}^2.
\end{equation}
\end{lem}
\noindent{\it Proof of the lemma \ref{lem52}.} Simply write the integral as a sum of low frequencies $\{\xi\in\mathbb R^3;\;|\xi|<1\}$ and high frequencies $\{\xi\in\mathbb R^3;\;|\xi|>1\}$.
\section{Well-posedness of $(NSE)$ in $H^1_{a,\sigma}(\mathbb{R}^3)$}
In the following theorem, we characterize the existence and uniqueness of the solution.
\begin{theorem}\label{theo2} Let $a>0$ and $\sigma\geq1$.
Let $u^0\in \left(H^1_{a,\sigma}(\mathbb{R}^3)\right)^3$ such that ${\rm div}\, u^0=0$. Then, there exist a unique $T^* \in(0,+\infty]$ and a unique solution $u\in{\mathcal C}([0,T^*),H^1_{a,\sigma}(\mathbb R^3))$ to $(NSE)$ system such that $u\notin{\mathcal C}([0,T^*],H^1_{a,\sigma}(\mathbb R^3))$. Moreover, if $T^*$ is finite, then
\begin{eqnarray}\label{enq1}
\lim\sup_{t\nearrow T^*}\|u(t)\|_{\dot{H}^1_{a,\sigma}}=\infty.
\end{eqnarray}
\end{theorem}
\noindent{\it Proof theorem \ref{theo2}.} The integral form of Navier-Stokes equations is
$$
u=e^{\nu t\Delta}u^0+B(u,u),
$$
where
$$
B(u,v)=-\int_0^t e^{\nu (t-\tau)\Delta}\mathbb{P}({\rm div}\,(u\otimes v))
$$
and
$$
\mathcal{F}(\mathbb{P}f)(\xi)=\widehat{f}(\xi)-\frac{(f(\xi).\xi)}{|\xi|^2}\xi.
$$
The main idea of this proof is to apply the theorem of the fixed point to the operator
$$u\longmapsto e^{\nu t\Delta}u^0+B(u,u).$$
Taking account of lemma \ref{lem52}, it suffices to estimate the nonlinear part respectively in the spaces $L^2(\mathbb R^3)$ and $\dot H^1_{a,\sigma}(\mathbb R^3)$.\\\\
$\bullet$ To estimate $B(u,v)$ in $\dot{H}^1_{a,\sigma}(\mathbb{R}^3)$, we have
\begin{eqnarray*}
\|B(u,v)\|_{\dot{H}^1_{a,\sigma}}(t)&\leq&\int_0^t\| e^{\nu (t-\tau)\Delta}\mathbb{P}({\rm div}\,(u\otimes v))\|_{\dot{H}^1_{a,\sigma}}\\ &\leq&\int_0^t\left(\int_{\xi}|\xi|^3  e^{-2\nu (t-\tau)|\xi|^2}|\xi| e^{2a|\xi|^{\frac{1}{\sigma}}}|\widehat{u\ast v}(\tau,\xi)|^2\right)^{\frac{1}{2}}\\&\leq&\int_0^t\frac{1}{\nu^{\frac{3}{4}}(t-\tau)^{\frac{3}{4}}}\|u\otimes v\|_{\dot{H}^{\frac{1}{2}}_{a,\sigma}},
\end{eqnarray*}
with
\begin{eqnarray*}
\|u\otimes v\|_{\dot{H}^{\frac{1}{2}}_{a,\sigma}}^2(\tau)&=&\int_{\xi} |\xi|e^{2a|\xi|^{\frac{1}{\sigma}}}|\widehat{u\ast v}(\tau,\xi)|^2\\&\leq&\int_{\xi}|\xi|\left(\int_{\eta} e^{a|\xi|^{\frac{1}{\sigma}}}|\widehat{u}(\tau,\xi-\eta)||\widehat{v}(\tau,\eta)|\right)^2.
\end{eqnarray*}
Using the inequality $e^{a|\xi|^{\frac{1}{\sigma}}}\leq e^{a|\xi-\eta|^{\frac{1}{\sigma}}}e^{a|\eta|^{\frac{1}{\sigma}}}$, we get
\begin{eqnarray*}
\|u\otimes v\|_{\dot{H}^{\frac{1}{2}}_{a,\sigma}}^2(\tau)&\leq&\int_{\xi}|\xi|\left(\int_{\eta} e^{a|\xi-\eta|^{\frac{1}{\sigma}}}|\widehat{u}(\tau,\xi-\eta)|e^{a|\eta|^{\frac{1}{\sigma}}}|\widehat{v}(\tau,\eta)|\right)^2\\
&\leq&\|UV\|_{\dot{H}^{\frac{1}{2}}}\\
&\leq&c\|U\|_{\dot{H}^1}\|V\|_{\dot{H}^1}\;\;({\rm by\,lemma\,\ref{lem1}})\\
&\leq&c\|u\|_{\dot{H}^1_{a,\sigma}}\|v\|_{\dot{H}^1_{a,\sigma}},
\end{eqnarray*}
where
$\widehat{U}=e^{a|\xi|^{\frac{1}{\sigma}}}|\widehat{u}(\tau,\xi)|$ and $\widehat{V}=e^{a|\xi|^{\frac{1}{\sigma}}}|\widehat{v}(\tau,\xi)|$.\\
Therefore,
\begin{eqnarray}\label{enq2}
\|B(u,v)\|_{\dot{H}^1_{a,\sigma}}(t)\leq c\nu^{-\frac{3}{4}}T^{\frac{1}{4}}\|u\|_{L^{\infty}_T(\dot{H}^1_{a,\sigma})}\|v\|_{L^{\infty}_T(\dot{H}^1_{a,\sigma})}.
\end{eqnarray}
\noindent$\bullet$ To estimate $B(u,v)$ in $L^2(\mathbb{R}^3)$, we have
\begin{eqnarray*}
\|B(u,v)\|_{L^2}(t)&\leq&\int_0^t\| e^{\nu (t-\tau)\Delta}\mathbb{P}({\rm div}\,(u\otimes v))\|_{L^2}\\ &\leq&\int_0^t\left(\int_{\xi}e^{-2\nu (t-\tau)|\xi|^2}|\xi|^2 |\widehat{u\ast v}(\tau,\xi)|^2\right)^{\frac{1}{2}}\\
&\leq&\int_0^t\frac{1}{\nu^{\frac{1}{4}}(t-\tau)^{\frac{1}{4}}}\|u\otimes v\|_{\dot H^{\frac{1}{2}}},\\
&\leq&\int_0^t\frac{1}{\nu^{\frac{1}{4}}(t-\tau)^{\frac{1}{4}}}\|u\|_{\dot H^1}\|v\|_{\dot H^1}\;\;({\rm by\,lemma\,\ref{lem1}}),\\
&\leq&\int_0^t\frac{1}{\nu^{\frac{1}{4}}(t-\tau)^{\frac{1}{4}}}\|u\|_{\dot H^1_{a,\sigma}}\|v\|_{\dot H^1_{a,\sigma}},\\
&\leq&\nu^{-\frac{1}{4}}t^{\frac{3}{4}}\|u\|_{\dot H^1_{a,\sigma}}\|v\|_{\dot H^1_{a,\sigma}}.
\end{eqnarray*}
Therefore,
\begin{eqnarray}\label{enq22}
\|B(u,v)\|_{L^2}(t)\leq c\nu^{-\frac{1}{4}}T^{\frac{3}{4}}\|u\|_{L^{\infty}_T(\dot{H}^1_{a,\sigma})}\|v\|_{L^{\infty}_T(\dot{H}^1_{a,\sigma})}.
\end{eqnarray}
Now, we recall the following fixed point theorem.
\begin{lem}\label{lem3}(See \cite{MC1})
Let $X$ be an abstract Banach space equipped with the norm $\|\,\|$ and let $B:X\times X\rightarrow X$ be a bilinear operator, such that for any $x_1,x_2\in X$,
$$
\|B(x_1,x_2)\|\leq c_0\|x_1\|\|x_2\|.
$$
Then, for any $y\in X$ such that
$$
4c_0\|y\|<1
$$
the equation
$$
x=y+B(x,x)
$$
has a solution $x\in X$. In particular, this solution satisfies
$$
\|x\|\leq2\|y\|
$$
and it is the only one such that
$$
\|x\|<\frac{1}{2c_0}.
$$
\end{lem}
Combining the inequalities (\ref{enq2})-(\ref{enq22}) and lemma \ref{lem3}, and choosing a time $T$ small enough, we guarantee the existence and uniqueness of the solution in the space ${\mathcal C}([0,T],H^1_{a,\sigma}(\mathbb R^3))$.\\\\

\noindent$\bullet$ Now, we want to prove the inequality (\ref{enq2}) by contradiction. Let $u$ be maximal solution to $(NSE)$ with $T^*<\infty$. It suffices to observe that $u$ is bounded, then we can extend the solution to an interval $[0,T_1]$, where $T_1>T^*$.\\
As the function $u$ is assumed to be bounded, then there exists $M\geq0$ such that
$$
\|u(t)\|_{H^1_{a,\sigma}(\mathbb{R}^3)}\leq M,\,\,\forall t\in[0,T^*].
$$
We deduce that, $\forall \varepsilon>0 \,\,\exists\alpha/\,\,\forall t,t'<T^*$ such that $T^*-t<\alpha$ and $T^*-t'<\alpha$, we have
$$
\|u(t)-u(t')\|_{\dot{H}^1_{a,\sigma}(\mathbb{R}^3)}<\varepsilon.
$$
Then $(u(t))$ is a Cauchy sequence in ${H^1_{a,\sigma}(\mathbb{R}^3)}$ at the time $T^*$. The space $H^1_{a,\sigma}(\mathbb{R}^3)$ is complete, then there exists $u^*\in H^1_{a,\sigma}(\mathbb{R}^3)$ such that
$$
\lim_{t\nearrow T^*}u(t)=u^*\;\;{\rm in}\; H^1_{a,\sigma}(\mathbb{R}^3).
$$
Consider the following system
$$
\left\{
  \begin{array}{lll}
     \partial_t v
 -\nu\Delta v+ v.\nabla v  & =&-\nabla p' \hbox{ in }  \mathbb R^+\times \mathbb R^3\\
    v(0) &=& u^*   \hbox{ in }\mathbb R^3.
  \end{array}
\right.\leqno{(NSE^*)}
$$
Applying the first step to $(NSE^*)$ system, then there exist $T_1>0$ and a unique solution $v\in{\mathcal C}([0,T_1],H^1_{a,\sigma}(\mathbb R^3))$ to $(NSE^*)$. So
$$
w(t)=\left\{
  \begin{array}{lll}
     u(t) &\hbox{ if }& t\in [0,T^*]\\
    v(t-T^*) &\hbox{ if }& t\in [T^*,T^*+T_1]
  \end{array}
\right.
$$
as a unique solution extends $u$, which is absurd. Taking into account the energy inequality (\ref{ee1}) we can deduce equation (\ref{enq1}). Therefore theorem \ref{theo2} is proved.
\section{Proof of the main result}
\subsection{Fundamental inequalities}
In this section, we give some blow up results which are necessary to prove our main result.
\begin{theorem}\label{theo3}
Let $u^0\in \left(H^1_{a,\sigma}(\mathbb{R}^3)\right)^3$ such that ${\rm div}\, u^0=0$ and $u\in{\mathcal C}([0,T^*),H^1_{a,\sigma}(\mathbb R^3))$ the maximal solution to $(NSE)$ system given by theorem \ref{theo1}. If $T^*<\infty$, then
\begin{eqnarray}\label{enq3}
\int_t^{T^*}\|\mathcal{F}(e^{\frac{a}{\sigma}|D|^{\frac{1}{\sigma}}}u(\tau))\|^2_{ L^1}d\tau=\infty,\,\,\forall t\in[0,T^*)
\end{eqnarray}
and
\begin{eqnarray}\label{enq4}
\frac{(\frac{\nu}{2})^{\frac{1}{2}}}{(T^*-t)^{\frac{1}{2}}}\leq\|\mathcal{F}(e^{\frac{a}{\sigma}|D|^{\frac{1}{\sigma}}}u)(t)\|_{ L^1},\,\,\forall t\in[0,T^*).
\end{eqnarray}
\end{theorem}
\noindent{\it Proof of theorem \ref{theo3}.}
$\bullet$ First, we wish to prove the inequality (\ref{enq3}).\\
Taking the scalar product in $\dot{H}^1_{a,\sigma}(\mathbb R^3)$, we obtain
\begin{eqnarray*}
\frac{1}{2}\partial_t \|u(t)\|_{\dot{H}^1_{a,\sigma}}^2+\nu\|\nabla u(t)\|_{\dot{H}^1_{a,\sigma}}^2&\leq&|\langle e^{a|D|^{\frac{1}{\sigma}}}( u.\nabla u)/e^{a|D|^{\frac{1}{\sigma}}}u(t) \rangle_{\dot{H}^1}|\\&\leq&|\langle e^{a|D|^{\frac{1}{\sigma}}}{\rm div}\, (u\otimes u)/e^{a|D|^{\frac{1}{\sigma}}}u(t) \rangle_{\dot{H}^1}|\\&\leq&|\langle e^{a|D|^{\frac{1}{\sigma}}} u\otimes u /\nabla e^{a|D|^{\frac{1}{\sigma}}}u(t)\rangle_{\dot{H}^1}|\\&\leq&\|u\otimes u\|_{\dot{H}^1_{a,\sigma}}\|\nabla u\|_{\dot{H}^1_{a,\sigma}}.
\end{eqnarray*}
Using lemma \ref{lem4}, we have
\begin{eqnarray*}
\frac{1}{2}\partial_t \|u(t)\|_{\dot{H}^1_{a,\sigma}}^2+\nu\|\nabla u(t)\|_{\dot{H}^1_{a,\sigma}}^2&\leq& 32 \|\mathcal{F}(e^{\frac{a}{\sigma}|D|^{\frac{1}{\sigma}}}u(t))\|_{ L^1}\|u(t)\|_{\dot{H}^1_{a,\sigma}}\|\nabla u(t)\|_{\dot{H}^1_{a,\sigma}}.
\end{eqnarray*}
Inequality $xy\leq \frac{x^2}{2}+\frac{y^2}{2}$ gives us
\begin{eqnarray*}
\frac{1}{2}\partial_t \|u(t)\|_{\dot{H}^1_{a,\sigma}}^2+\nu\|\nabla u(t)\|_{\dot{H}^1_{a,\sigma}}^2&\leq&c\nu^{-1}\|\mathcal{F}(e^{\frac{a}{\sigma}|D|^{\frac{1}{\sigma}}}u(t))\|_{ L^1}^2\|u(t)\|_{\dot{H}^1_{a,\sigma}}^2+\frac{\nu}{2}\|\nabla u(t)\|_{\dot{H}^1_{a,\sigma}}^2.
\end{eqnarray*}
Then
\begin{eqnarray*}
\partial_t \|u(t)\|_{\dot{H}^1_{a,\sigma}}^2+\nu\|\nabla u(t)\|_{\dot{H}^1_{a,\sigma}}^2&\leq&c\nu^{-1}\|\mathcal{F}(e^{\frac{a}{\sigma}|D|^{\frac{1}{\sigma}}}u(t))\|_{ L^1}^2\|u(t)\|_{\dot{H}^1_{a,\sigma}}^2.
\end{eqnarray*}
Gronwall lemma implies that, for every $0\leq t\leq T<T^*$,
$$
\|u(T)\|_{\dot{H}^1_{a,\sigma}}^2\leq\|u(t)\|_{\dot{H}^1_{a,\sigma}}^2 e^{c\nu^{-1}\int_t^T\|\mathcal{F}(e^{\frac{a}{\sigma}|D|^{\frac{1}{\sigma}}}u(t))\|_{ L^1}^2}.
$$
The fact that $\lim\sup_{t\nearrow T^*}\|u(t)\|_{\dot{H}^1_{a,\sigma}}=\infty$ implies that
\begin{eqnarray*}
\int_t^{T^*}\|\mathcal{F}(e^{\frac{a}{\sigma}|D|^{\frac{1}{\sigma}}}u(\tau))\|^2_{ L^1}=\infty,\,\,\forall t\in[0,T^*).
\end{eqnarray*}
$\bullet$ Now, we want to prove inequality (\ref{enq4}).\\
Returning to $(NSE)$ system and taking the Fourier transform in the first equation, we get
$$
\frac{1}{2}\partial_t|\widehat{u}(t,\xi)|^2+\nu |\xi|^2|\widehat{u}(t,\xi)|^2+Re\left(\mathcal{F}(u.\nabla u)(t,\xi).\widehat{u}(t,-\xi)\right)=0.
$$
But for any $\varepsilon>o$,
$$
\frac{1}{2}\partial_t|\widehat{u}(t,\xi)|^2=\frac{1}{2}\partial_t(|\widehat{u}(t,\xi)|^2+\varepsilon)=\sqrt{|\widehat{u}(t,\xi)|^2
+\varepsilon}.\partial_t\sqrt{|\widehat{u}(t,\xi)|^2+\varepsilon}
$$
and
$$
\partial_t\sqrt{|\widehat{u}(t,\xi)|^2+\varepsilon}+\nu|\xi|^2\frac{|\widehat{u}(t,\xi)|^2}{\sqrt{|\widehat{u}(t,\xi)
|^2+\varepsilon}}+\frac{\mathcal{F}(u.\nabla u)(t,\xi).\widehat{u}(t,-\xi)}{\sqrt{|\widehat{u}(t,\xi)|^2+\varepsilon}}=0.
$$
Then, we have
$$
\partial_t\sqrt{|\widehat{u}(t,\xi)|^2+\varepsilon}+\nu|\xi|^2\frac{|\widehat{u}(t,\xi)|^2}{\sqrt{|\widehat{u}(t,\xi)|^2+\varepsilon}}\leq
|\mathcal{F}(u.\nabla u)(t,\xi)|.
$$
Integrating on $[t,T]\subset[0,T^*)$, we obtain
$$
\sqrt{|\widehat{u}(T,\xi)|^2+\varepsilon}+\nu|\xi|^2\int_t^T\frac{|\widehat{u}(\tau,\xi)|^2}{\sqrt{|\widehat{u}(\tau,\xi)|^2+\varepsilon}}
\leq\sqrt{|\widehat{u}(t,\xi)|^2+\varepsilon}+\int_t^T|\mathcal{F}(u.\nabla u)(\tau,\xi)|.
$$
Taking $\varepsilon\rightarrow 0$, we get
\begin{eqnarray*}
|\widehat{u}(T,\xi)|+\nu|\xi|^2\int_t^T|\widehat{u}(\tau,\xi)|&\leq&|\widehat{u}(t,\xi)|+\int_t^T|\mathcal{F}(u.\nabla u)(\tau,\xi)|\\&\leq&|\widehat{u}(t,\xi)|+\int_t^T|\widehat{u}|\ast_{\xi}|\widehat{\nabla u}|(\tau,\xi) \\&\leq&|\widehat{u}(t,\xi)|+\int_t^T\int_{\eta}|\widehat{u}(\tau,\xi-\eta)||\widehat{\nabla u}(\tau,\eta)|.
\end{eqnarray*}
Multiplying the late equation by $e^{\frac{a}{\sigma}|\xi|^{\frac{1}{\sigma}}}$ and using the inequality
$$
e^{\frac{a}{\sigma}|\xi|^{\frac{1}{\sigma}}}\leq e^{\frac{a}{\sigma}|\xi-\eta|^{\frac{1}{\sigma}}}e^{\frac{a}{\sigma}|\eta|^{\frac{1}{\sigma}}},
$$
we obtain
\begin{eqnarray*}
e^{\frac{a}{\sigma}|\xi|^{\frac{1}{\sigma}}}|\widehat{u}(T,\xi)|+\nu|\xi|^2\int_t^T e^{\frac{a}{\sigma}|\xi|^{\frac{1}{\sigma}}}|\widehat{u}(\tau,\xi)|\\ \leq e^{\frac{a}{\sigma}|\xi|^{\frac{1}{\sigma}}}|\widehat{u}(t,\xi)|+\int_t^T\int_{\eta}e^{\frac{a}{\sigma}|\xi-\eta|^{\frac{1}{\sigma}}}|\widehat{u}
(\tau,\xi-\eta)|e^{\frac{a}{\sigma}|\eta|^{\frac{1}{\sigma}}}|\widehat{\nabla u}(\tau,\eta)|.
\end{eqnarray*}
Integrating over  $\xi\in \mathbb R^3$ and using Young's inequality, we get
\begin{eqnarray*}
\|\mathcal{F}(e^{\frac{a}{\sigma}|D|^{\frac{1}{\sigma}}}u)(T)\|_{ L^1}+\nu\int_t^T\|\mathcal{F}(e^{\frac{a}{\sigma}|D|^{\frac{1}{\sigma}}}\Delta u)\|_{ L^1}\\
\leq\|\mathcal{F}(e^{\frac{a}{\sigma}|D|^{\frac{1}{\sigma}}}u)(t)\|_{ L^1}+\int_t^T\|\mathcal{F}(e^{\frac{a}{\sigma}|D|^{\frac{1}{\sigma}}}u)\|_{ L^1}\|\mathcal{F}(e^{\frac{a}{\sigma}|D|^{\frac{1}{\sigma}}}\nabla u)\|_{ L^1}.
\end{eqnarray*}
Cauchy-Schwarz inequality yields
\begin{eqnarray*}
\|\mathcal{F}(e^{\frac{a}{\sigma}|D|^{\frac{1}{\sigma}}}u)(T)\|_{ L^1}+\nu\int_t^T\|\mathcal{F}(e^{\frac{a}{\sigma}|D|^{\frac{1}{\sigma}}}\Delta u)\|_{ L^1}\\
\leq\|\mathcal{F}(e^{\frac{a}{\sigma}|D|^{\frac{1}{\sigma}}}u)(t)\|_{ L^1}+\int_t^T\|\mathcal{F}(e^{\frac{a}{\sigma}|D|^{\frac{1}{\sigma}}}u)\|^{\frac{3}{2}}_{ L^1}\|\mathcal{F}(e^{\frac{a}{\sigma}|D|^{\frac{1}{\sigma}}}\Delta u)\|^{\frac{1}{2}}_{ L^1}.
\end{eqnarray*}
Inequality  $xy\leq \frac{x^2}{2}+\frac{y^2}{2}$ gives
\begin{eqnarray*}
\|\mathcal{F}(e^{\frac{a}{\sigma}|D|^{\frac{1}{\sigma}}}u)(T)\|_{ L^1}+\frac{\nu}{2}\int_t^T\|\mathcal{F}(e^{\frac{a}{\sigma}|D|^{\frac{1}{\sigma}}}\Delta u)\|_{L^1}\leq\|\mathcal{F}(e^{\frac{a}{\sigma}|D|^{\frac{1}{\sigma}}}u)(t)\|_{L^1}+\frac{1}{\nu}\int_t^T\|\mathcal{F}(e^{\frac{a}{\sigma}|D|^{\frac{1}{\sigma}
}}u)\|^3_{ L^1}.
\end{eqnarray*}
By the Gronwall lemma, for $0\leq t\leq T<T^*$, we obtain
$$
\|\mathcal{F}(e^{\frac{a}{\sigma}|D|^{\frac{1}{\sigma}}}u)(T)\|_{ L^1}\leq \|\mathcal{F}(e^{\frac{a}{\sigma}|D|^{\frac{1}{\sigma}}}u)(t)\|_{ L^1}e^{\nu^{-1}
\int_t^T\|\mathcal{F}(e^{\frac{a}{\sigma}|D|^{\frac{1}{\sigma}}}u)(\tau)\|^2_{ L^1}}
$$
or
$$
\|\mathcal{F}(e^{\frac{a}{\sigma}|D|^{\frac{1}{\sigma}}}u)(T)\|^2_{L^1}e^{-2\nu^{-1}\int_t^T\|\mathcal{F}(e^{\frac{a}{\sigma}|D|^{\frac{1}{\sigma}}}u)
(\tau)\|^2_{ L^1}}\leq  \|\mathcal{F}(e^{\frac{a}{\sigma}|D|^{\frac{1}{\sigma}}}u)(t)\|^2_{ L^1}.
$$
Integrating over $[t_0,T]\subset[0,T^*)$, we get
$$
1-e^{-2\nu^{-1}\int_t^T\|\mathcal{F}(e^{\frac{a}{\sigma}|D|^{\frac{1}{\sigma}}}u)(\tau)\|^2_{L^1}}\leq2\nu^{-1}\|\mathcal{F}(e^{\frac{a}{\sigma}
|D|^{\frac{1}{\sigma}}}u)(t_0)\|^2_{ L^1}(T-t_0).
$$
By inequality (\ref{enq3}), if $T\rightarrow T^*$, we have
$$
1\leq 2\nu^{-1}\|\mathcal{F}(e^{\frac{a}{\sigma}|D|^{\frac{1}{\sigma}}}u)(t_0)\|^2_{ L^1}(T^*-t_0).
$$
Then we can deduce inequality (\ref{enq4}), and theorem \ref{theo3} is proved.
\subsection{Exponential type explosion} Let $a>0$ and $\sigma>1$. Let $u\in{\mathcal C}([0,T^*), H^1_{a,\sigma}(\mathbb R^3))$ be a maximal solution to $(NSE)$ system given by theorem \ref{theo2} and suppose that $T^*<\infty$. We want to prove inequality (\ref{ineq99}).\\

\noindent$\bullet$ Firstly, we prove the following result
\begin{equation}\label{eq999}
\frac{\frac{\nu}{2}}{T^*-t}\leq\|\widehat{u}(t)\|^2_{ L^1},\,\,\forall t\in[0,T^*).
\end{equation}
For this, using Cauchy-Schwarz inequality, we get
\begin{eqnarray*}
\|\mathcal{F}(e^{\frac{a}{\sigma}|D|^{\frac{1}{\sigma}}}u)(t)\|_{L^1}&=&\int_{\xi}e^{\frac{a}{\sigma}|\xi|^{\frac{1}{\sigma}}}
|\widehat{u}(\xi)|\\&\leq& c_{a,\sigma}\|u\|_{\dot{H}^1_{\frac{a}{\sqrt{\sigma}},\sigma}}
\end{eqnarray*}
where
\begin{eqnarray*}
c^2_{a,\sigma}=\int_{\mathbb R^3}\frac{1}{|\xi|^2}e^{-2a(\frac{1}{\sqrt{\sigma}}-\frac{1}{\sigma})|\xi|^{\frac{1}{\sigma}}}=4\pi\sigma \big(2a(\frac{1}{\sqrt{\sigma}}-\frac{1}{\sigma})\big)^{\sigma-2}\Gamma(\sigma)<\infty.
\end{eqnarray*}
Since $\frac{a}{\sqrt{\sigma}}<a$, then $\dot{H}^1_{a,\sigma}(\mathbb R^3))\hookrightarrow \dot{H}^1_{\frac{a}{\sqrt{\sigma}},\sigma}(\mathbb R^3))$.\\
Let's take  $a'=\frac{a}{\sqrt{\sigma}}$, we ensure that $u\in{\mathcal C}([0,T^*),\dot{H}^1_{a',\sigma}(\mathbb R^3))$ and
$$
\lim\sup_{t\nearrow T^*}\|u(t)\|_{\dot{H}^1_{a',\sigma}}=\infty.
$$
Using similar technique, we get
$$
\int_t^{T^*}\|\mathcal{F}(e^{\frac{a'}{\sqrt{\sigma}}|D|^{\frac{1}{\sigma}}}u(\tau))\|^2_{ L^1}=\infty,\,\,\forall t\in[0,T^*)
$$
and
$$
\frac{\frac{\nu}{2}}{T^*-t}\leq\|\mathcal{F}(e^{\frac{a'}{\sqrt{\sigma}}|D|^{\frac{1}{\sigma}}}u)(t)\|^2_{ L^1},\,\,\forall t\in[0,T^*).
$$
This implies
$$
\frac{\frac{\nu}{2}}{T^*-t}\leq\|\mathcal{F}(e^{\frac{a}{(\sqrt{\sigma})^2}|D|^{\frac{1}{\sigma}}}u)(t)\|^2_{ L^1},\,\,\forall t\in[0,T^*).
$$
By induction, we can deduce that for every $n\in\mathbb{N}$,
$$
\frac{\frac{\nu}{2}}{T^*-t}\leq\|\mathcal{F}(e^{\frac{a}{(\sqrt{\sigma})^n}|D|^{\frac{1}{\sigma}}}u)(t)\|^2_{ L^1},\,\,\forall t\in[0,T^*).
$$
Using dominated convergence theorem, we get inequality (\ref{ineq99}) by letting $n\rightarrow +\infty$.\\

\noindent$\bullet$ Secondly, we prove the exponential type explosion. Using lemma \ref{lem2} and the energy estimate
$$\|u(t)\|_{L^2}^2+2\nu\int_0^t\|\nabla u(z)\|_{L^2}^2dz\leq \|u^0\|_{L^2}^2,$$
we can write, for $k\in\mathbb N$ such that $1+\frac{k}{2\sigma}\geq2$ (i.e $k\geq2\sigma$),
$$\frac{\sqrt{\frac{\nu}{2}}}{\sqrt{T^*-t}}\leq\|\widehat{u}(t)\|_{ L^1}\leq M(2)\|u^0\|_{L^2}^{1-\frac{3}{2(1+\frac{k}{2\sigma})}}\|u(t)\|_{\dot H^{1+\frac{k}{2\sigma}}}^{\frac{3}{2(1+\frac{k}{2\sigma})}},\,\,\forall t\in[0,T^*).
$$
Then
$$\Big(\frac{\sqrt{\frac{\nu}{2}}}{\sqrt{T^*-t}}\Big)^{2(\frac{1+\frac{k}{2\sigma}}{3})}
(M(2))^{-2(\frac{1+\frac{k}{2\sigma}}{3})}\leq\|u^0\|_{L^2}^{2(\frac{1+\frac{k}{2\sigma}}{3})-1}\|u(t)\|_{\dot H^{1+\frac{k}{2\sigma}}}
$$
or
$$\Big(\frac{\sqrt{\frac{\nu}{2}}}{\sqrt{T^*-t}}\Big)^{2(\frac{1+\frac{k}{2\sigma}}{3})}
(M(2))^{-2(\frac{1+\frac{k}{2\sigma}}{3})}\|u^0\|_{L^2}^{1-2(\frac{1+\frac{k}{2\sigma}}{3})}\leq\|u(t)\|_{\dot H^{1+\frac{k}{2\sigma}}}.
$$
We obtain
$$
\frac{C_1}{(T^*-t)^{2/3}}\Big(\frac{C_2}{(T^*-t)^{1/3\sigma}}\Big)^k\leq \|u(t)\|_{\dot H^{1+\frac{k}{2\sigma}}}^2
$$
with
$$ C_1=\left(\frac{\nu}{2}(M(2))^{-2}\|u^0\|_{L^2}\right)^{\frac{2}{3}},\,\,C_2=\left(\frac{\nu}{2}(M(2))^{-2}\|u^0\|_{L^2}^{2}\right)^{\frac{1}{3\sigma}}.$$
Then
$$\frac{1}{k!}\frac{C_1}{(T^*-t)^{2/3}}\Big(\frac{2aC_2}{(T^*-t)^{1/3\sigma}}\Big)^k\leq \int_\xi\frac{(2a)^k}{k!}|\xi|^{\frac{k}{\sigma}}|\xi|^2|\widehat{u}(t,\xi)|^2.$$
Summing over the set $\{k\in\mathbb N;\;\;k\geq 2\sigma\}$, we get
$$\frac{C_1}{(T^*-t)^{1/3}}\Big(e^{\frac{2aC_2}{(T^*-t)^{1/3\sigma}}}-\sum_{0\leq k\leq 2\sigma}\frac{(\frac{2aC_2}{(T^*-t)^{1/3\sigma}})^k}{k!}\Big)\leq \int_\xi\Big(e^{2a|\xi|^{\frac{1}{\sigma}}}-\sum_{0\leq k\leq 2\sigma}\frac{(2a|\xi|^{\frac{1}{\sigma}})^k}{k!}\Big)|\xi|^2|\widehat{u}(t,\xi)|^2$$
$$\leq\int_\xi e^{2a|\xi|^{\frac{1}{\sigma}}}|\xi|^2|\widehat{u}(t,\xi)|^2.$$
Now, put $$h(z)=\displaystyle\frac{e^z-\displaystyle\sum_{k=0}^{2\sigma_0}\frac{z^k}{k!}}{z^{2\sigma_0+1}e^\frac{z}{2}},\;\;z>0,$$
with $2\sigma_0$ is the integer part of $2\sigma$.\\\\
The function $h$ satisfies the following properties\\
$\bullet$ $h$ is continuous on $(0,\infty)$,\\
$\bullet$ $h(z)>0$ for all $z>0$,\\
$\bullet$ $\displaystyle \lim_{z\rightarrow\infty}h(z)=\infty$,\\
$\bullet$ $\displaystyle \lim_{z\rightarrow0^+}h(z)=\frac{1}{(2\sigma_0+1)!}$.\\
Then there is $B=B(\sigma_0)>0$ such that $h(z)\geq B$ for all $z>0$. Therefore\\
$$\frac{c_1}{(T^*-t)^{\frac{1}{3}+\frac{2\sigma_0+1}{3\sigma}}}\exp\Big[\frac{aC_2}{(T^*-t)^{1/3\sigma}}\Big]\leq \|u(t)\|_{H^1_{a,\sigma}},$$
with $c_1=BC_1(2aC_2)^{2\sigma_0+1}.$ Then, we can deduce inequality (\ref{ineq99}) and theorem \ref{theo1} is proved.

\end{document}